\documentclass[12pt,twoside,final]{amsart}
\usepackage[psamsfonts]{amssymb}
\usepackage{times,a4wide}
\usepackage{esint}
\usepackage[sc]{mathpazo}
\usepackage{xcolor}
\usepackage{url}
\usepackage{bbm}
\usepackage{orcidlink}

\usepackage{amsrefs}
\usepackage{hyperref}

\theoremstyle{plain}
\newtheorem*{theorem*}{Theorem}

\newtheorem{theorem}{Theorem}[section]
\newtheorem{proposition}[theorem]{Proposition}
\newtheorem*{maintheorem*}{Main Theorem}
\newtheorem*{proposition*}{Proposition}
\newtheorem{corollary}[theorem]{Corollary}
\newtheorem*{corollary*}{Corollary}
\newtheorem{lemma}[theorem]{Lemma}
\newtheorem*{lemma*}{Lemma}

\theoremstyle{definition}

\newtheorem*{remark*}{Remark}
\newtheorem*{remarks*}{Remarks}
\newtheorem*{conjecture*}{Conjecture}
\newtheorem*{Chernoff}{Chernoff's bounds}
\theoremstyle{definition}
\newenvironment{poliabstract}[1]
  {\renewcommand{\abstractname}{#1}\begin{abstract}}
  {\end{abstract}}

\allowdisplaybreaks

\usepackage[active]{srcltx}
\usepackage{tikz}
\usepackage{pgfplots}
\usepackage{float}
\RequirePackage{pgfplots}
\usetikzlibrary{shadows}
\pgfplotsset{compat=1.13} % XAVIER : 1.16
\newcommand{\C}{\mathbb{C}}

\newcommand{\N}{\mathbb{N}}
\newcommand{\R}{\mathbb{R}}
\newcommand{\E}{\mathbb{E}}

\renewcommand{\P}{\mathbb{P}}

\newcommand{\Cov}{\operatorname{Cov}}

\renewcommand{\arg}{\operatorname{arg}}

\title[Separation properties of a hybrid process] 
{Separation properties of a hybrid point process with determinantal radii and uniform arguments} 

\author[G. Lamberti]{Giuseppe Lamberti
\orcidlink{0009-0009-0503-0421}}
\address{
Center for Mathematical Analysis, Geometry and Dynamical Systems \&
Departamento de Matem\'atica\\
Instituto Superior T\'ecnico\\
Universidade de Lisboa\\
Av. Rovisco Pais\\
1049-001 Lisboa, Portugal.}
\email{giuseppe.lamberti@tecnico.ulisboa.pt}

\author[X. Massaneda]{Xavier Massaneda
\orcidlink{0000-0002-7552-9190}}

\address{Departament de Matem\`atiques i Inform\`atica,
Universitat  de Barcelona, Gran Via 585, 08007-Bar\-ce\-lo\-na, Catalonia}
\email{xavier.massaneda@ub.edu}

\thanks{Giuseppe Lamberti was funded by FCT/Portugal and the Recovery and Resilience Plan (PRR) through projects UID/04459/2025 and UID/PRR/04459/2025, and by the project 2023.17881.ICDT (SHADE). The first author is partially supported by the project UBGRS 2.0 (ANR-20-SFRI-0001) and by the grant ``INdAM-GNAMPA Project'', CUP E53C25002010001, ``Transferring Harmonic Analysis between Discrete Structures and Manifolds''. Second author partially supported by the Generalitat de Catalunya (grant 2021 SGR 00087) and the Spanish Ministerio de Ciencia e Innovaci\'on (project PID2024-160033NB-I00).}

\date{\today}


\begin{document}

%\selectlanguage{english}
\begin{poliabstract}{Abstract} 
We recently characterized the almost sure separation of determinantal point processes $\Lambda_\phi$ in the plane associated with Fock spaces $\mathcal F_\phi$ with doubling weights $\phi$. We also showed that, as expected, a more restrictive condition is required to characterize the almost surely separated Poisson processes with the same first intensity as $\Lambda_\phi$. To gain further insight into this difference in behavior, we focus on radial weights $\phi$ and introduce a hybrid process $\Lambda_\phi^M=\{r_j e^{i\theta_j}\}_{j}$, where the moduli $r_j$ are taken from $\Lambda_\phi$, while the arguments $\theta_j$ are chosen independently and uniformly in $[0,2\pi)$, independently of the radii. Our main result is that the hybrid process is almost surely separated under exactly the same condition as the corresponding Poisson process. Thus, despite preserving the joint radial distribution of the determinantal process, independent angular randomization destroys its improved separation properties.
\end{poliabstract}

\maketitle

\section{Introduction}\label{sec:intro}

In this paper we characterize the almost sure separation of a family of point processes in the complex plane.  Recall that a sequence $\Lambda=\{\lambda_j\}_{j}\subset\C$ is (Euclidean) \emph{separated} if
\[
 \inf_{j\neq k} |\lambda_k-\lambda_j|>0.
\]
Separation plays an important role in various problems in analysis, particularly in the characterization of interpolating and sampling sequences for spaces of analytic functions.

Throughout the paper, we consider only simple point processes, that is, random locally finite configurations of distinct points in $\mathbb C$. Equivalently, such a process can be viewed as a random counting measure
\[
 \tau_\Lambda=\sum_{\lambda\in \Lambda}\delta_\lambda ,
\]
where $\Lambda$ is a locally finite subset of $\C$. Its distribution (or law) is determined by the joint distributions of random counting variables
\[
N(B)=\#(\Lambda\cap B)=\int_B d\tau_\Lambda,\qquad \text{$B\subset\C$ compact}.
\]
For general background on point processes and random measures, see \cite{Da-Ve}.

In a recent paper we characterized the doubling subharmonic weights $\phi$ for which the determinantal process $\Lambda_\phi$ naturally associated with the generalized Fock space $\mathcal F_\phi$ is almost surely separated \cite{LaMa}*{Theorem 2.1}. We also showed that, as expected, a more restrictive condition is required for a Poisson process $\Lambda_\phi^P$ with the same first intensity as $\Lambda_\phi$ to be almost surely separated.

In order to state the results, recall that a subharmonic weight is called \emph{doubling} when its Laplacian $\nu=\Delta\phi$ is a doubling measure, that is, there exists $C>0$ such that 
\[
 \nu(D(z,2r))\leq C\, \nu(D(z,r)) ,\quad\text{for all}\ z\in\C,\quad r>0.
\]
We will always assume that $\nu(\C)=\infty$; otherwise the process is finite.

The results are expressed in terms of the scale function $\rho(z)$ (or $\rho_\phi(z)$ if we want to stress the dependence on $\phi$), which is the unique positive radius such that
\begin{equation}\label{eq:defrho}
 \nu(D(z,\rho(z)))=1.
\end{equation}
The function $\rho^{-2}$ can then be seen as a regularized version of $\Delta\phi$, as described by M. Christ in \cite{Christ}, and is comparable to the first intensity of the process (see Section~\ref{sec-det}).
%The scale function $\rho$ describes the local density of the process: its first intensity is comparable to $\rho^{-2}$.
%Actually, according to Theorem 14 in \cite{MMOC}, given a subharmonic doubling weight $\phi$, there exists $\psi$ subharmonic, doubling and smooth such that $|\psi-\phi|$ is uniformly bounded and $\Delta\psi\simeq \rho_\phi^{-2}\simeq \rho_\psi^{-2}$. 

In \cite{LaMa}*{Theorem 2.1} we proved that, for the determinantal process, 
\[
\Lambda_\phi\ \text{a.s. separated}\quad \Longleftrightarrow \quad \int_{\C}\frac{dm}{\rho^6}<+\infty.
\]
We also showed \cite{LaMa}*{Theorem 2.5} that, for the Poisson process $\Lambda_\phi^P$ having the same first intensity as $\Lambda_\phi$,
\[
\Lambda_\phi^P\ \text{a.s. separated}\quad \Longleftrightarrow \quad \int_{\C}\frac{dm}{\rho^4}<+\infty.
\]

The different exponents in these two conditions reflect the stronger local repulsion exhibited by determinantal processes. This raises a natural question: to what extent is this repulsion carried by the moduli of the points, and to what extent by their arguments?

\iffalse
In this paper we investigate the separation of the hybrid process 
\[
\Lambda_\phi^M=\{r_j e^{i\theta_j}\}_{j}, 
\]
where the moduli $r_j$ are taken from the determinantal process associated to a radial weight $\phi$ (for precise definitions see next section), and the arguments $\theta_j$ are i.i.d. uniform variables in $[0,2\pi)$, and independent of the radii. 
\fi

To better understand the situation, we introduce a hybrid process that preserves the joint distribution of the moduli of the determinantal process $\Lambda_\phi$, while independently randomizing their arguments. More precisely, for a radial doubling subharmonic weight $\phi$, we define
\[
 \Lambda_\phi^M=\{r_j e^{i\theta_j}\}_j,
\]
where $\{r_j\}_j $ have the same distribution as $\{|\lambda| : \lambda\in\Lambda_\phi\}$ and the $\theta_j$ are i.i.d. uniform random variables on $[0,2\pi)$, independent of $\{r_j\}_j $.

By rotational invariance, $\Lambda_\phi^M$ has the same first intensity as $\Lambda_\phi$ and, consequently, as the Poisson process $\Lambda_\phi^P$. Thus, the three processes have identical first intensities but different correlation structures.

Our main result, Theorem~\ref{thm:main}, shows that, as far as almost sure separation is concerned, the hybrid process behaves like the Poisson process. Thus, the radial structure inherited from the determinantal process does not improve the separation threshold: after independent angular randomization, the hybrid process satisfies exactly the same separation criterion as the Poisson process.

\begin{theorem}\label{thm:main} Let $\phi$ be a radial doubling subharmonic function, and let $\Lambda_\phi^M$ be the associated hybrid point process. Then
\[
 \P\bigl(\Lambda_\phi^M\ \text{is separated}\bigr)=
 \begin{cases}
  1\quad \text{if $\displaystyle{\int_{\C}}\dfrac{dm(z)}{\rho^4(z)}<+\infty$},\\
  0\quad\text{if $\displaystyle{\int_{\C}}\dfrac{dm(z)}{\rho^4(z)}=+\infty$}.
 \end{cases}
\]
\end{theorem}

Canonical examples of the weights considered here are $\phi_\alpha(z)=|z|^\alpha$, $\alpha>0$. The particular case $\alpha=2$ corresponds to the classical Fock space.  Since $\Delta|z|^\alpha=\alpha^2 |z|^{\alpha-2}$ and $\rho(z)\simeq |z|^{1-\alpha/2}$ for $|z|\geq 1$, one sees that 
$\int_{\C} \rho^{-4} dm<\infty$ if and only if $\int_1^\infty r^{2\alpha-3} dr<\infty$; hence:

\begin{corollary}\label{cor:alpha} Let $\Lambda_\alpha^M$ be the hybrid point process associated with the weight $\phi_\alpha(z)=|z|^\alpha$, $\alpha>0$. Then, 
\[
 \P\bigl(\Lambda_\alpha^M\ \text{is separated}\bigr)=
 \begin{cases}
  1\quad \text{if $\alpha<1$},\\
  0\quad\text{if $\alpha\geq 1$}.
 \end{cases}
\]
\end{corollary}

In contrast, the determinantal process $\Lambda_\alpha$ is almost surely separated if and only if $\alpha<4/3$ \cite{LaMa}. Thus, for $1\leq\alpha<4/3$, the determinantal process is almost surely separated, while the corresponding hybrid and Poisson processes are almost surely not separated.

The same arguments as in \cite{LaMa}*{Subsection 3.4} show that when $\Lambda_\phi^M$ is separated (actually, whenever $\rho(x)\to\infty$ as $x\to\infty$) its upper density
\[
 D_+(\Lambda_\phi^M):=\limsup_{r\to\infty}\sup_{z\in\C} \frac{\#\bigl(\Lambda_\phi^M\cap D(z,r)\bigr)}{\pi r^2}
\]
is zero. Consequently, by the classical characterization of interpolating sequences for Fock spaces given by K. Seip and R. Wallst\'en \cites{Seip92a,Seip92b}, $\Lambda_\phi^M$ is almost surely an interpolating sequence for $\mathcal F_{\gamma|z|^2}$ for every $\gamma>0$.

\begin{corollary}\label{cor:interpolation} Let $\phi$ be a radial doubling subharmonic weight and let $\Lambda_\phi^M$ be its associated hybrid process. Then, for any $\gamma>0$,
\[
 \P\bigl(\Lambda_\phi^M\ \text{is $\mathcal F_{\gamma|z|^2}$-interpolating}\bigr)=
 \begin{cases}
   1\quad \text{if $\displaystyle{\int_{\C}}\dfrac{dm(z)}{\rho^4(z)}<+\infty$},\\
  0\quad\text{if $\displaystyle{\int_{\C}}\dfrac{dm(z)}{\rho^4(z)}=+\infty$}.
 \end{cases}
\]
\end{corollary}

The plan of the paper is as follows.
Section~\ref{sec-det} contains the detailed description of the hybrid probabilistic model and its basic properties. In Section~\ref{sec-red}, we reduce the proof of the main theorem to a discretized separation problem (Proposition~\ref{prop:main-partition}), under the assumption that $\rho(x)\to\infty$ as $x\to\infty$. Section~\ref{sec-proof:main-prop} develops the probabilistic estimates needed for the more delicate divergent case of Theorem~\ref{thm:main} (or of Proposition~\ref{prop:main-partition}), including negative-association arguments and estimates for collision probabilities.

Finally, the Appendix deals with the elementary remaining case in which $\rho(x)$ does not tend to infinity as $x\to\infty$, completing the proof of the main theorem.

A final word about notation: throughout the paper, $A\lesssim B$ means that $A\leq CB$ for some constant $C>0$ independent of the relevant variables. We write $A\simeq B$ if both $A\lesssim B$ and $B\lesssim A$ hold.

\section{The hybrid probabilistic model}\label{sec-det} 

Consider the Fock space of entire functions associated with $\phi$ subharmonic and doubling
\[
\mathcal F_{\phi}=\bigl\{f\in H(\C) : \|f\|_\phi^2:=\int_{\C}|f(z)|^2 e^{-2\phi(z)}\, \frac{dm(z)}{\rho_\phi^2 (z)}<+\infty\bigr\}.
\]
Equipped with the scalar product
\[
 \langle f,g \rangle_{\phi}=\int_{\C} f(z)\, \overline{g(z)}\, e^{-2\phi(z)}\, \frac{dm(z)}{\rho_\phi^2 (z)},\qquad f,g\in \mathcal F_{\phi},
\]
the space $\mathcal F_{\phi}$ is a reproducing kernel Hilbert space whose kernel is denoted by $K_\phi$ (see for instance \cite{MMOC}). Specifically, this means that $K_\phi$ is holomorphic in $z$, anti-holomorphic in $\zeta$, and $f(z)=\langle f,K_\phi(z,\cdot)\rangle_{\phi}$ for any $f\in \mathcal F_\phi$.

Denote 
\begin{equation}\label{eq:measure}
 d\mu_\phi(z)=e^{-2\phi(z)}\, \frac{dm(z)}{\rho_\phi^2 (z)}.
\end{equation}

According to a theorem of Macchi and Soshnikov (\cite{soshnikov00}*{Theorem 3}, see also \cite{HKPV}*{Lemma~4.5.1}), there exists a determinantal point process $\Lambda_\phi$ associated with $K_\phi$ and $\mu_\phi$. This means that for any collection of disjoint sets $B_1,\dots, B_n\subset\C$ one has
\[
 \mathbb E\left[\prod_{k=1}^n \#(B_k\cap\Lambda)\right]=
 \int_{B_1}\cdots \int_{B_n}\det\bigl(K_\phi(z_i,z_j)\bigr)_{1\leq i,j\leq n} d\mu_\phi(z_1)\cdots d\mu_\phi(z_n).
\]
The integrand is usually called the $n^{\text{th}}$-\textit{correlation function} of $\Lambda_\phi$. For the particular case $n=1$, one has
\[
 \mathbb E[N(B)]=\int_{B} K_\phi(z,z)\, d\mu_\phi(z),
\]
and the measure $K_\phi(z,z)\, d\mu_\phi(z)$ is called the \textit{first intensity} (or average distribution) of the process.

Since, by Lemma 21 in \cite{MMOC}, there exists $C>0$ such that
\begin{equation}\label{eq:kernel-diagonal}
C^{-1} e^{2\phi(z)}\leq K_\phi(z,z)\leq C\, e^{2\phi(z)}\  \quad\text{for all}\ z\in\C,
\end{equation} 
the first intensity of $\Lambda_\phi$ is comparable to the regularization of the Laplacian of $\phi$: 
\begin{equation}\label{eq:1st-intensity}
 K_\phi (z,z)\, d\mu_\phi(z)\simeq \frac{dm(z)}{\rho^2(z)}.
\end{equation}
In particular, given any Borel set $B\subset \C$, the expected value of the counting variable $N_\phi(B)=\#(\Lambda_\phi\cap B)$ is
\begin{equation}\label{eq:EB}
 \mathbb E[N_\phi(B)]=\int_B K_\phi (z,z)\, d\mu_\phi(z)\simeq \int_B\frac{dm(z)}{\rho^2(z)}.
\end{equation}

According to \cite{HKPV}*{Theorem 4.7.1}, given a determinantal point process $\Lambda$ with kernel $K(z,w)$ and background radial finite measure $\mu$, the radii $\{|\lambda|:\lambda\in \Lambda\}$ have the same law as the set $\{r_j\}_{j=0}^\infty$ described next.

Since the background measure is radial, the kernel has the form
\[
 K(z,w)=\sum_{j=0}^\infty a_j^2 (z\bar w)^j,
\]
where $e_j(z)=a_j z^j$ are the normalized eigenfunctions for $K$. Denote by $\varphi(|z|)$ the density of the measure $\mu$. Let $\{r_j\}_j$ be independent positive random variables such that $r_j^2$ has the distribution $Q_j$ given by the density function
\[
 f_j(t)=\pi a_j^2 t^j \varphi(\sqrt t),\qquad t>0.
\]
In our case, by \eqref{eq:measure}, the density is
\[
 \varphi(r)= \frac{e^{- 2\phi(r)}}{ \rho^2(r)},\qquad r>0.
\]
Since
\begin{align*}
 \|z^j\|_{\phi}^2 &=\int_{\C} |z|^{2j} e^{-2\phi(z)}\frac{dm(z)}{\rho^2(z)}=\pi \int_0^\infty r^{2j} e^{-2\phi(r)}\frac {2r\, dr}{\rho^2(r)}
 =\pi \int_0^\infty t^{j} \varphi(\sqrt t)\, dt,
\end{align*}
one has
\begin{equation}\label{eq:aj}
 \pi a_j^2=\left(\int_0^\infty t^{j} \varphi(\sqrt t)\, dt\right)^{-1},
\end{equation}
and therefore
\begin{equation}\label{eq:density-rk}
f_j(t)=\frac{t^j\, \varphi(\sqrt t)}{\int_0^\infty s^{j} \varphi(\sqrt s)\, ds} ,\qquad t>0.
\end{equation}

The hybrid point process we consider in this paper is $\Lambda_\phi^M=\{\lambda_j^M\}_{j\geq 0}=\{r_j e^{i\theta_j}\}_{j\geq 0}$, where the radii are independent, $r_j$ has density $f_j$, and $\theta_j$ are i.i.d uniform in $[0,2\pi)$, and independent of the radii. 

%%%%%%%%%%%%%%%%%%%%%%%%%%%%%%%%%%%%%%%%%%%%%%%%%%%%%%%%%%%%%%%%%%%%%%%%%%%%%%%

\section{Discretization and reduction to Proposition~\ref{prop:main-partition}}\label{sec-red}

The main theorem (Theorem~\ref{thm:main}) can be reduced to a simpler, discretized version
by considering a standard partition of the plane and estimating the ``collision'' probability that two or more points fall into the same cell.

For $n\geq 1$ and $k=1,\dots, n$, consider the annuli
\[
 I_n=\{z\in\C : n-1\leq |z| <n\}
\]
and the angular cells
\[
T_{n,k}:= \Big\{z \in \C:\ n-1\leq |z| <n, \ \frac{\arg(z)}{2\pi} \in \big[\frac{k-1}{n}, \frac{ k}{n}\big)\ (\text{mod 1})\Big\}. 
\]
Denote by $N_n=\#(\Lambda_{\phi}^M\cap I_n)$
and by $X_{n,k}=\#(\Lambda_{\phi}^M\cap T_{n,k})$
the corresponding counting random variables.

Notice that the random variable $N_n$ can  be written as the sum 
\begin{equation}\label{eq:Nn-Bernoulli}
 N_n=\sum_{j\geq 0} \zeta_j^{(n)},
\end{equation}
where $\zeta_j^{(n)}\sim\text{Bernoulli}(p_j^{(n)})$ are independent and each $\zeta_j^{(n)}$ indicates whether 
$\lambda_j^{M}$ falls in $I_n$ or not.
Using the densities of $r_j^2$ given in \eqref{eq:density-rk} we see that:
\begin{align}\label{eq:pkn}
 p_j^{(n)}&=\mathbb P\bigl(\zeta_j^{(n)}=1\bigr)=
 \mathbb P\bigl(r_j\in [(n-1)^2,n^2)\bigr)=
 \frac{\int_{(n-1)^2}^{n^2} t^j\, \varphi(\sqrt t)\, dt}{\int_{0}^\infty t^j\, \varphi(\sqrt t)\, dt}\\
 &=\frac{\int_{n-1}^n r^{2j+1}\, \varphi(r)\, dr}{\int_{0}^\infty r^{2j+1}\, \varphi(r)\, dr}
 =\frac{\int_{n-1}^n r^{2j+1}\, e^{-2\phi(r)}\, \frac{dr}{\rho^2(r)}}{\int_{0}^\infty r^{2j+1}\,  e^{-2\phi(r)}\, \frac{dr}{\rho^2(r)}}. \notag
\end{align}

Observe that the $\{p_j^{(n)}\}_j$ are precisely the eigenvalues of the restriction operator
\begin{equation*}
 Tf(z)=\int_{I_{n}} f(\zeta) K_{\phi}(z,\zeta) d\mu_{\phi}(\zeta),\qquad f\in\mathcal F_{\phi},
\end{equation*}
because the normalized monomials $e_j(z)=a_j z^j$ are complete and orthogonal in all $I_n$ and, by \eqref{eq:aj},
\begin{align*}
 Te_j(z)&=\int_{I_{n}} e_j(\zeta) \sum_{m=0}^\infty e_m(z) \overline{e_m(\zeta)} \, d\mu_\phi(\zeta)=e_j(z) \int_{I_{n}} |e_j(\zeta)|^2 \, d\mu_\phi(\zeta)\\
 &=e_j(z) 2\pi a_j^2\int_{n-1}^n r^{2j+1} e^{-2\phi(r)}\frac {dr}{\rho^2(r)}=p_j^{(n)}\, e_j(z).
\end{align*}
Thus, $N_n$ can also be thought of as a sum of independent Bernoulli variables whose parameters are the eigenvalues of the restriction operator $T_n$ (see \cite{HKPV}*{Theorem 4.5.3}).

Before examining the distributions of $N_n$ and $X_{n,k}$, we recall here some well-known properties of the radius $\rho(z)$ defined in \eqref{eq:defrho}. 

\begin{lemma}[\cite{MMOC}*{Section 2.1}]\label{lem:rho} Let $\phi$ be a doubling subharmonic weight and let $\rho(z)$ denote the radius defined in \eqref{eq:defrho}. Then
\begin{itemize}
 \item [(a)] $|\rho(z)-\rho(\zeta)|\leq |z-\zeta|$, for any $ z,\zeta\in \C$.
 
 \item [(b)] There exist $C>0$ and $\beta\in (0,1)$ such that
\begin{equation}\label{eq:beta}
 \rho(z)\leq C |z|^\beta\qquad |z|\geq 1.
\end{equation}

\item [(c)] For every $r>1$ there exists $C_r>1$ such that for $\zeta\in D(z,r\rho(z))$
\[
\frac 1{C_r}\leq\frac{\rho(z)}{\rho(\zeta)}\leq C_r.
\]
\end{itemize}
\end{lemma}

When studying the separation of $\Lambda_\phi^M$, the only interesting case occurs when 
 \begin{equation}\label{eq:rho-to-inft}
  \liminf_{x\to +\infty}\rho(x)=+\infty,
 \end{equation}
which we assume henceforth (see Appendix or Lemma 3.8 in \cite{LaMa}). This corresponds to the cases where the first intensity is asymptotically smaller than the area measure, and determining whether $\Lambda_\phi^M$ is a.s. separated becomes more delicate. It is also worth pointing out that if \eqref{eq:rho-to-inft} fails, then $\int_{\C} \rho^{-4} dm=\infty$, as shown in \cite{LaMa}*{Section 3.2, formula (13)}.

Under \eqref{eq:rho-to-inft}, by \eqref{eq:1st-intensity} and Lemma~\ref{lem:rho}(c), it is clear that
\begin{equation}\label{eq:ENn}
 \mu_n:=\mathbb E[N_n]=\sum_{j\geq 0} p_j^{(n)}\simeq \int_{I_n}\frac{dm(z)}{\rho^2(z)}\simeq \frac{|I_n|}{\rho^2(n)}\simeq \frac{n}{\rho^2(n)}.
\end{equation}

We can then rewrite the critical integral in the statement of Theorem~\ref{thm:main} as
\begin{equation}\label{eq:int-mun}
 \int_{\C} \frac{dm(z)}{\rho^4(z)}=\sum_{n=1}^\infty \int_{I_n} \frac{dm(z)}{\rho^4(z)}\simeq\sum_{n=1}^\infty\frac {n}{\rho^4(n)}
 \simeq\sum_{n=1}^\infty\frac {\mu_n^2}{n}.
\end{equation}
Furthermore, by the uniformity in the distributions of the angles in $\lambda_k^M$, the distribution of the random variables $X_{n,k}$, $k=1,\dots,n$ does not depend on $k$. In particular
\[
 \mathbb E[X_{n,k}]=\frac{\mathbb E[N_n]}n\simeq \frac{1}{\rho^2(n)}.
\]

In this situation, we also have good estimates of the probabilities $p_j^{(n)}$.

\begin{lemma}\label{lem:pnk}
Let $\phi$ be a radial doubling weight and assume that $\rho$ satisfies \eqref{eq:rho-to-inft}. Let $p_j^{(n)}$ be the probabilities defined in \eqref{eq:pkn}.

\begin{itemize}
 \item [(a)] There exists $C_0>0$ such that for all $n\geq 1$
 \begin{equation}\label{eq:bound-pnk}
  \sup_{j\geq 0} p_j^{(n)}\leq\frac {C_0}{\rho(n)}.
 \end{equation}
 
 \item [(b)] $\displaystyle{\lim_{n\to\infty}\frac{\sup_{j\geq 0} p_j^{(n)}}{\mu_n}}=0$.
 \end{itemize}
\end{lemma}

\begin{proof}
Set
\[
    H_j(r)=\frac {r^{2j+1} e^{-2\phi(r)}}{\rho^2(r)}. 
\]

(a) Fix $n$ and $j$, and take $s_{n,j}\in[n-1,n]$ such that
\[
 H_j(s_{n,j})=\max_{r\in[n-1,n]} H_j(r),
\]
 so that
 \begin{equation}\label{eq:max}
  \int_{n-1}^n  H_j(r)\,  dr \leq  H_j(s_{n,j}).
 \end{equation}
 By Lemma 19 in \cite{MMOC}, for every $R>0$ there exists $C_0=C_0(R)>0$ such that for any $f\in H(\C)$ and any $z\in\C$
 \[
  |f(z)|^2 e^{-2\phi(z)}\leq C_0\int_{D(z,R\rho(z))} |f(\zeta)|^2 e^{-2\phi(\zeta)}\, \frac{dm(\zeta)}{\rho^2(\zeta)}.
 \]
 Applying this to $R=1$, $f(z)=z^j$ and the point $z=s_{n,j}$, we see that
 \[
  s_{n,j}^{2j} e^{-2\phi(s_{n,j})}\leq C_0\int_{D(s_{n,j},\rho(s_{n,j}))} |\zeta|^{2j} e^{-2\phi(\zeta)}\, \frac{dm(\zeta)}{\rho^2(\zeta)}.
 \]
 Since, by Lemma~\ref{lem:rho}(b), $\rho(s_{n,j})\ll s_{n,j}$ for $n$ large enough, we can replace the disk $D(s_{n,j},\rho(s_{n,j}))$ by the bigger angular sector 
 \[
  Q(s_{n,j},\rho(s_{n,j}))=\bigl\{re^{i\theta} : |r-s_{n,j}|<\rho(s_{n,j})\ ;\ |\theta|<\arcsin(\frac{\rho(s_{n,j})}{s_{n,j}})\bigr\}
 \]
 and integrate in polar coordinates; it follows that
 \begin{align*}
   s_{n,j}^{2j}\, e^{-2\phi(s_{n,j})}&\leq C_0\int_{Q\bigl(s_{n,j},\rho(s_{n,j})\bigr)} |\zeta|^{2j} e^{-2\phi(\zeta)} \frac{dm(\zeta)}{\rho^2(\zeta)}\\
&   =2C_0\bigl(\arcsin(\frac{\rho(s_{n,j})}{s_{n,j}})\bigr) 
   \int_{s_{n,j}-\rho(s_{n,j})}^{s_{n,j}+\rho(s_{n,j})} H_j(r)\, dr.
 \end{align*}
Again by Lemma~\ref{lem:rho}(b), 
\[
 \arcsin(\frac{\rho(s_{n,j})}{s_{n,j}})\simeq \frac{\rho(s_{n,j})}{s_{n,j}}\qquad n\to\infty,
\]
so we deduce from the inequality above that
\[
\int_{s_{n,j}-\rho(s_{n,j})}^{s_{n,j}+\rho(s_{n,j})} H_j(r)\, dr
\gtrsim  \rho(s_{n,j})\, H_j(s_{n,j}).
\]
This, together with \eqref{eq:max}, shows finally that
\begin{align*}
 p_j^{(n)}&=\frac{\int_{n-1}^n H_j(r)\, dr}{\int_{0}^\infty H_j(r)\, dr}\leq 
 \frac{\int_{n-1}^n H_j(r)\, dr}{\int_{s_{n,j}-\rho(s_{n,j})}^{s_{n,j}+\rho(s_{n,j})} H_j(r)\, dr}\lesssim \frac 1{\rho(s_{n,j})}.
\end{align*}
Since $\rho(n)\to\infty$ and $|s_{n,j}-n|\leq 1$, we can use Lemma~\ref{lem:rho}(c) to deduce that $\rho(s_{n,j})\simeq \rho(n)$, with implicit constants independent of $n$, so that (a) follows.

(b) is an immediate consequence of (a) and Lemma~\ref{lem:rho}(b):
\[
 \frac{\sup_{j\geq 0} p_j^{(n)}}{\mu_n}\lesssim \frac{1/\rho(n)}{n/\rho^2(n)}=\frac{\rho(n)}n\longrightarrow 0,\qquad n\to\infty.
\]
\end{proof}

Let us recall now the standard scheme, dating back at least to the proof of Theorem 2 in \cite{Coch} (see also \cite{LaMa}*{Section 3.3}), to reduce Theorem~\ref{thm:main} to the following discretized version.

For every $n\geq 1$ consider the collision event that at least one of the cells in $I_n$ has two points or more
\begin{equation}\label{eq:En}
E_n := \bigl\{\exists\, k = 1,\dots,n : X_{n,k} \geq 2\bigr\} = \Bigl\{\max_{1 \leq k \leq n} X_{n,k} \geq 2\Bigr\}.
\end{equation}
Observe that for $m\neq n$ the events $E_n$ and $E_m$ are not independent: a given point cannot lie in both $I_n$ and $I_m$, so the events are dependent. In fact, as we shall see, these events are negatively correlated (see Proposition~\ref{prop:na}).

\begin{proposition}\label{prop:main-partition}
Let $\phi$ be a doubling radial weight satisfying \eqref{eq:rho-to-inft}
and let $X_{n,k}$ be as above. Then
\[
 \P\bigl(E_n\ \text{infinitely often}\bigr)=
 \begin{cases}
  0\quad  \text{if $\displaystyle{\sum_{n=1}^\infty\frac {\mu_n^2}{n}}<+\infty$}\\
  1 \quad  \text{if $\displaystyle{\sum_{n=1}^\infty\frac {\mu_n^2}{n}}=+\infty$.}
 \end{cases}
\]
\end{proposition}

We now deduce Theorem \ref{thm:main} from Proposition \ref{prop:main-partition}.

\begin{proof}[Proof of Theorem \ref{thm:main}]

The convergent case $\sum_{n=1}^\infty \mu_n^2/n <+\infty$ is easy: from the first Borel-Cantelli lemma we deduce that almost surely $X_{n,k}\leq 1$ for all but at most finitely many $n,k$. 

Technically, this still does not imply that $\Lambda_\phi^M$ is separated, since points in adjacent boxes could be arbitrarily close to one another. But the arguments of the proof of Proposition~\ref{prop:main-partition} (see Section \ref{sec-proof:main-prop}) can be applied analogously to shifted regions, both in the radius and the arguments, which allows then to deduce the a.s. separation. Let us see briefly how this can be done.

For $a,b\in \{0,1/2\}$ consider the sectors
\begin{align*}
 T_{n,k}^{a,b}=\Big\{z\in\C : n-1+a\leq |z|<n+a,\, \frac{\arg(z)}{2\pi} \in \big[\frac{k-1+b}{n}, \frac{ k+b}{n}\big)\, (\text{mod $1$})\Big\}, 
\end{align*}
and the corresponding random variables $ X_{n,k}^{a,b}=\#(\Lambda_{\phi}^M\cap  T_{n,k}^{a,b})$. It will be clear from the proof that Proposition~\ref{prop:main-partition} remains valid if the events $E_n$ are replaced by $E_n^{a,b}=\bigl\{\exists\, k = 1,\dots,n : X_{n,k}^{a,b} \geq 2\bigr\}$, i.e. $\P(E_n^{a,b}\ \text{infinitely often})=0$. Therefore, with probability one $X_{n,k}^{a,b}\leq 1$ for any combination of $a,b\in\{0,1/2\}$ and for all but at most finitely many $n,k$. 

It is not difficult to see that if $|z_1|,|z_2|\geq 2$ and $|z_1-z_2|<1/10$, then there exist $(n,k)$ and $a,b$ such that $z_1,z_2\in T_{n,k}^{a,b}$. Thus, except for maybe a finite number of points, the sequence $\Lambda_\phi^M$ is separated (by the constant $\delta=1/10$).

For the case $\sum_{n=1}^\infty \mu_n^2/n =+\infty$, the arguments of the proof of Proposition~\ref{prop:main-partition} can be applied similarly to any grid of size $h\in(0,1]$: letting, for $n\geq 1$ and $k=1,\dots, n$,
\[
 T_{n,k}^h:= \Big\{z \in \C: (n-1)h\leq |z|<nh, \, \frac{\arg(z)}{2\pi} \in \big[\frac{k-1}{n}, \frac{ k}{n}\big)\, (\text{mod $1$})\Big\}, 
\]
and $X_{n,k}^{h}=\#(\Lambda_{\phi}^M\cap T_{n,k}^h)$, we see that the events
\[
 E_n^h=\{\exists k\in\{1,\dots,n\} : X_{n,k}^{h}\geq 2\}
\]
happen infinitely often with probability 1. 

Taking $h=1/l$, $l\geq 1$, and defining $\mathcal E_l=\{E_n^{1/l}\, \text{infinitely often}\}$ we obtain an event
\[
 \mathcal E=\cap_{l\geq 1} \mathcal E_l
\]
which has probability one and under which, for any $l\geq 1$, there are infinitely many pairs of points inside the same cell $T_{n,k}^{1/l}$. In particular,  $\Lambda_\phi^M$ is not separated.
\end{proof}

\section{Proof of Proposition~\ref{prop:main-partition}}\label{sec-proof:main-prop}
Consider the collision events $E_n$ defined in \eqref{eq:En}. From the first Borel-Cantelli lemma, it is clear that if $\sum_n \mathbb P(E_n)<\infty$ then $\P(E_n\, \text{infinitely often}) =0$.

On the other hand, since the events $E_n$ are not independent, we cannot apply the second Borel-Cantelli and we shall instead rely on the following generalized version. 

\begin{lemma}\label{lem:BC2}
 If for all $n\neq m$
 \[
  \P(E_n\cap E_m)\leq \P(E_n)\, \P(E_m)
 \]
 and $\sum_n \P(E_n)=+\infty$, then
 \[
  \P( E_n\, \text{infinitely often})=1.
 \]
\end{lemma}
Notice that the extra hypothesis in this lemma can be rephrased in terms of the indicator functions $\mathbbm{1}_{E_n}$ as
\[
 \Cov(\mathbbm{1}_{E_n}, \mathbbm{1}_{E_m})\leq 0\qquad n\neq m.
\]

\begin{proof}
 For $N\leq M$, set $S_{N,M}=\sum\limits_{n=N}^M \mathbbm{1}_{E_n}$ and denote $s_{N,M}=\E [S_{N,M}]=\sum_{n=N}^M \P(E_n)$.
 By hypothesis
 \begin{align*}
  \E[S_{N,M}^2]&=\E\Bigl[\sum\limits_{n=N}^M \mathbbm{1}_{E_n}+ \sum_{
  \substack{n,m=N \\ m\neq n}} \mathbbm{1}_{E_n} \mathbbm{1}_{E_m}\Bigr]\\
  &=s_{N,M}+ \sum_{\substack{n,m=N \\ m\neq n}} \Cov(\mathbbm{1}_{E_n}, \mathbbm{1}_{E_m})
  + \sum_{\substack{n,m=N \\ m\neq n}} \E[\mathbbm{1}_{E_n}] \E[\mathbbm{1}_{E_m}]\\
  &\leq s_{N,M}+s_{N,M}^2.
 \end{align*}
 By the Cauchy-Schwarz inequality,
 \[
  s_{N,M}^2=\bigl(\E [S_{N,M}]\bigr)^2\leq \E[S_{N,M}^2]\, \P(S_{N,M}\geq 1),
 \]
 hence
 \[
  \P(S_{N,M}\geq 1)\geq \frac{s_{N,M}^2}{\E[S_{N,M}^2]}\geq \frac{s_{N,M}^2}{s_{N,M}+s_{N,M}^2}=\frac{s_{N,M}}{1+s_{N,M}}.
 \]
By the divergence of $\sum_n \P(E_n)$,
\[
 \lim_{M\to\infty} \P(S_{N,M}\geq 1)= 1,
\]
thus $\P(\cup_{n\geq N} E_n)=1$ for all $N$ and the conclusion follows.
\end{proof}

We shall see next that the family of indicator random variables $(\mathbbm{1}_{E_n})_{n\in\N}$ is negatively associated (see definition below), which implies that the hypothesis of Lemma~\ref{lem:BC2} is satisfied (Proposition~\ref{prop:na}).

\subsection{Negatively associated variables}
A family $(Y_j)_{j\in J}$ of random variables is called \emph{negatively associated} (NA) if,
for all disjoint finite sets of indices $A,B\subset J$ and all bounded functions
$f$ and $g$ that are nondecreasing in each coordinate,
\[
  \Cov\bigl(f((Y_i)_{i\in A}),\,g((Y_j)_{j\in B})\bigr)\le 0.
\]
Here $f$ is nondecreasing in each coordinate if $f(x)\le f(y)$
whenever $x_i\le y_i$ for every $i$.
This notion was introduced by Joag-Dev and Proschan
\cite{JoagDevProschan}, who also proved the following:
indicators of pairwise disjoint events form an NA family;
the union of mutually independent NA families is NA; and
nondecreasing functions applied to pairwise disjoint subfamilies of
an NA family are again NA. The first of these facts has a short
direct proof, which we include for completeness.

\begin{lemma}\label{lem:onehot}
Let $\lambda$ be a random point and let $\{T_n\}_n$ be a finite or
countable family of pairwise disjoint measurable sets. Then the
indicator variables $(\mathbbm{1}_{\{\lambda\in T_n\}})_n$ are
negatively associated.
\end{lemma}

\begin{proof}
Take disjoint finite index sets \(A,B\), and let \(f,g\) be bounded
nondecreasing functions of \(X_A=(\mathbbm{1}_{\{\lambda\in T_n\}})_{n\in A}\)
and \(X_B=(\mathbbm{1}_{\{\lambda\in T_m\}})_{m\in B}\). Since \(\lambda\) lies in
at most one set $T_n$, \(X_A\) is either \(0\) or a unit vector \(e_n\),
\(n\in A\), that is, a vector whose $n$-th coordinate is $1$ and all its other coordinates are $0$. Therefore, letting $\mathbf{0}$ denote the vector with all entries $0$,
\[
 f(X_A)=f(\mathbf{0})+\sum_{n\in A}\mathbbm{1}_{\{\lambda\in T_n\}}\,(f(e_n)-f(\mathbf{0})),
\]
and similarly for \(g(X_B)\). Set
\[
 a=\sum_{n\in A}\P(\lambda\in T_n)\,(f(e_n)-f(\mathbf{0})),
 \qquad
 b=\sum_{m\in B}\P(\lambda\in T_m)\,(g(e_m)-g(\mathbf{0})),
\]
and notice that $a,b\geq 0$, by the monotonicity of $f$ and $g$. Since \(\lambda\) cannot
belong to a cell from \(A\) and a cell from \(B\) at once, the cross
terms in \(f(X_A)g(X_B)\) vanish. Hence 
\[
\E[f(X_A)]=f(\mathbf{0})+a,\quad \E[g(X_B)]=g(\mathbf{0})+b, 
\]
and 
\[
\E[f(X_A)g(X_B)]=f(\mathbf{0})g(\mathbf{0})+g(\mathbf{0})a+f(\mathbf{0})b,
\]
so 
\[
\Cov(f(X_A),g(X_B))=\E[f(X_A)g(X_B)]-\E[f(X_A)]\E[g(X_B)]=-ab\leq0.
\]
\end{proof}

We now apply Lemma~\ref{lem:onehot} to show the negative correlation of the collision events, and hence verify the hypothesis of Lemma~\ref{lem:BC2}.

\begin{proposition}
\label{prop:na}
For the hybrid process, the events $(E_n)$ are
negatively correlated:
\[
 \Cov(\mathbbm{1}_{E_n},\mathbbm{1}_{E_m})\leq0\qquad(n\ne m).
\]
\end{proposition}

\begin{proof}
Set \(Y_{j,n,k}=\mathbbm{1}_{\{\lambda_j^M\in T_{n,k}\}}\), so that
$X_{n,k}=\sum_j Y_{j,n,k}$.

For each fixed \(j\),
Lemma~\ref{lem:onehot} shows that \((Y_{j,n,k})_{n,k}\) is NA. Since the
points \(\lambda_j^M\) are independent, these families are independent for
different \(j\), so the whole union \((Y_{j,n,k})_{j,n,k}\) is NA
\cite{JoagDevProschan}*{Property $P_7$}. 

The event \(E_n\) is nondecreasing in the
variables \(Y:=(Y_{j,n,k})_{j,k}\), meaning that the function 
\[
\mathbbm{1}_{E_n}(Y)=\mathbbm{1}_{\{\max\limits_k \sum_j Y_{j,n,k}\geq 2 \}}(Y)
\]
is a nondecreasing function in each coordinate separately.

Also, it is clear that for \(m\ne n\)  the function $\mathbbm{1}_{E_m}$
depends on a disjoint collection of variables $Y_{j,m,l}$. 

Since the definition of NA involves finitely many variables, we first restrict to the finite indices \(j\le N\):
the resulting events \(E_n^{(N)},E_m^{(N)}\) satisfy, by definition of negative association, the inequalities
\(\Cov(\mathbbm{1}_{E_n^{(N)}},\mathbbm{1}_{E_m^{(N)}})\leq 0\). As \(N\to\infty\), the
truncated counts and collision indicators increase to their full
counterparts, and bounded convergence gives the result.
\end{proof}

We have seen so far that

\begin{equation}\label{eq:S-finite}
 \mathbb P\bigl(E_n\ \text{infinitely often}\bigr)=
 \begin{cases}
  0 & \text{if $\sum_n \P(E_n)<+\infty$}\\
  1 & \text{if $\sum_n \P(E_n)=+\infty$}.
 \end{cases}
\end{equation}
Therefore, Proposition~\ref{prop:main-partition} will be proved as soon as we see that
$\sum_n \P(E_n)<+\infty$ if and only if $\sum_n \mu_n^2/n<+\infty$.

\subsection{Proof of Proposition~\ref{prop:main-partition}.
Case $\sum_{n=1}^\infty \mu_n^2/n <+\infty$}\label{subsec:finite-sum}
By \eqref{eq:S-finite}, it is enough to prove that $\sum_n \P(E_n)<+\infty$. Since the moduli $r_j$ and the arguments $\theta_j$ are independent, and $\theta_j$ are uniformly distributed, for all $j\geq 0$ and $k=1,\dots,n$,
\[
 \P\bigl(\lambda_j^M\in T_{n,k}\bigr)=\frac{p_j^{(n)}}{n},
\]
where $p_j^{(n)}=\P\bigl(\lambda_j^M\in I_{n}\bigr)$ is the probability given in \eqref{eq:pkn}.

If $E_n$ occurs, some cell $T_{n,k}$ contains two points $\lambda_i^M,\lambda_j^M$ with $i<j$, hence
\[
 E_n\subset\bigcup_{i<j}\ \bigcup_{k=1}^n\bigl\{\lambda_i^M\in T_{n,k},\ \lambda_j^M\in T_{n,k}\bigr\}.
\]
By the union bound and the independence of the points,
\begin{align*}
 \P(E_n)&\leq\sum_{i<j}\sum_{k=1}^n \P\bigl(\lambda_i^M\in T_{n,k}\bigr)\,\P\bigl(\lambda_j^M\in T_{n,k}\bigr)=\frac 1n\sum_{i<j}p_i^{(n)}p_j^{(n)}\\
 &\leq\frac 1{2n}\Bigl(\sum_{j\geq 0}p_j^{(n)}\Bigr)^2=\frac{\mu_n^2}{2n}.
\end{align*}
Therefore $\sum_{n=1}^\infty \P(E_n) \leq\frac 12\sum_{n=1}^\infty\mu_n^2/n<+\infty$.

\subsection{Proof of Proposition~\ref{prop:main-partition}\label{sec:3.2}.
Case $\sum_{n=1}^\infty \mu_n^2/n=+\infty$}\label{subsec:sum_n_div}

By \eqref{eq:S-finite}, it is enough to prove that $\sum_n \P(E_n)=+\infty$, for which we need a lower bound of $\P(E_n)$ (see Proposition~\ref{prop:collision} below). 

We first compute the probability of $E_n$ conditioned on the number of points in $I_n$, which is a classical birthday problem.  Denote
$q_{n,m}:=\P(E_n\mid N_n=m)$. By the independence and uniformity of the angles $\theta_j$,
\begin{equation}\label{eq:birthday}
 1-q_{n,m}=\P(X_{n,k}\leq 1\ \forall k \mid N_n=m)=
 \begin{cases}
  \frac{n(n-1)\cdots(n-m+1)}{n^m}\ &0\leq m\leq n,\\
  \ 0  &m>n.
 \end{cases}
\end{equation}
In particular
\(q_{n,0}=q_{n,1}=0\) and \(q_{n,2}=1/n\).

\begin{lemma}\label{lem:birthday}
There are absolute constants \(c_0,c_1>0\) such that, for all \(n\geq1\) and
integers \(m\geq0\),
\begin{equation}\label{eq:q}
 c_0\min\left\{1,\frac{m(m-1)}n\right\}
 \leq q_{n,m} \leq c_1 \min\left\{1,\frac{m(m-1)}n\right\}
\end{equation}
\end{lemma}

\begin{proof}
For the upper bound, label the $m$ points in $I_n$ as $\lambda_1,\dots,\lambda_m$, and for each pair $i<j$ let 
$A_{ij}=\{\lambda_i\ \text{and $\lambda_j$ fall in the same cell $T_{n,k}$}\}$. By definition, $\P(A_{ij})=1/n$ for every pair $i<j$. Then $E_n=\cup_{1\leq i<j\leq m} A_{ij}$ and
\[
 q_{n,m}\leq \sum_{i<j} \P(A_{ij})=\binom m2 \frac 1n.
\]

Conversely, for \(m\leq n\), equation \eqref{eq:birthday} gives
\[
 1- q_{n,m}
 =\prod_{i=0}^{m-1}\left(1-\frac in\right)
 \leq\exp\left(-\frac{m(m-1)}{2n}\right).
\]
Using the estimate \(1-e^{-x}\geq(1-e^{-1})\min\{1,x\}\), the stated inequality follows. 

The case \(m>n\) is immediate.
\end{proof}

Furthermore, we need the following classical bounds \cite{MiUp}.

\begin{Chernoff}[for the sum of independent Bernoulli random variables]
 Let $X$ be a sum of independent Bernoulli random variables and let $\mu=\E[X]$. Then:
 \begin{itemize}
  \item [(a)] $\P\bigl(X\geq (1+\delta)\mu\bigr)\leq e^{-\frac{\delta^2}{2+\delta}\, \mu}$, $\delta>0$.
  
  \item [(b)] $\P\bigl(X\leq (1-\delta)\mu\bigr)\leq e^{-\frac{\delta^2}{2}\, \mu}$, $\delta\in(0,1)$.
  
  \item [(c)] $\P\bigl(|X-\mu|\geq \delta\mu\bigr)\leq 2 e^{-\frac{\delta^2}{3}\, \mu}$, $\delta\in(0,1)$.
 \end{itemize}
\end{Chernoff}

These inequalities are usually stated for finite sums, but the proof shows that they hold as well for infinite sums, as long as $\mu=\E[X]<+\infty$.

We also recall the following result, due to Le Cam \cite{LeCa}. As with the Chernoff bounds, it is usually stated for finite sums, but it extends to infinite sums with $\mu<+\infty$ by a limiting argument.

\begin{theorem}[Le Cam's theorem] \label{thm:LeCam}
	Let $\{X_j\}_{j=1}^\infty$ be a sequence of independent Bernoulli random variables of parameter $p_j$ respectively. Suppose $\mu: = \sum_j p_j < \infty$ and define $S=\sum_j X_j$. Let $Y$ be a Poisson random variable of parameter $\mu$. Then
	\[
	\sum_{m=0}^\infty \left| \P(S=m)- \P(Y=m) \right| \leq 2\sum_{j=1}^\infty p_j^2
	\]
\end{theorem}

We can now state the estimate that will lead us to the conclusion of the proof.

\begin{proposition}\label{prop:collision}
There exist an absolute constant $c>0$ and $n_0\in\N$ such that
\begin{equation}\label{eq:lower}
 \P(E_n)\geq c\min\Bigl\{1,\frac{\mu_n^2}{n}\Bigr\},\qquad n\geq n_0.
\end{equation}
\end{proposition}

This finishes the proof of Proposition~\ref{prop:main-partition}, since $\sum_n\mu_n^2/n=+\infty$ yields $\sum_n \P(E_n)=\infty$.
%and the result follows from \eqref{eq:S-finite}.
 
\begin{proof}[Proof of Proposition~\ref{prop:collision}]
Denote $p_n^*=\sup_{j\geq 0} p_j^{(n)}$. By Lemma~\ref{lem:pnk}(b), there exists $n_0$ such that $p_n^*\leq (e^{-2}/8)\, \mu_n$ for all $n\geq n_0$. Fix such $n$; we distinguish two cases, depending on the size of $\mu_n$. 

\emph{{Case 1: $\mu_n\leq 2$.}}
Let $Y_n$ be a Poisson variable of parameter $\mu_n$. Since 
$\sum_j (p_j^{(n)})^2\leq p_n^* \,\mu_n$, 
Theorem~\ref{thm:LeCam} gives
\begin{align*}
 \P(N_n=2)&\geq\P(Y_n=2)-2\sum_{j\geq 0}\bigl(p_j^{(n)}\bigr)^2\geq e^{-\mu_n}\frac{\mu_n^2}2-2p_n^*\,\mu_n\\
 &\geq\Bigl(\frac{e^{-2}}2-\frac{e^{-2}}4\Bigr)\mu_n^2=\frac{e^{-2}}4\,\mu_n^2.
\end{align*}
Then, by \eqref{eq:birthday}, and since $\P(E_n\mid N_n=2)=q_{n,2}=1/n$,
\[
 \P(E_n)\geq \P\bigl(E_n\cap\{N_n=2\}\bigr)=\frac 1n\,\P(N_n=2)\geq\frac{e^{-2}}4\,\frac{\mu_n^2}n.
\]

\emph{{Case 2: $\mu_n>2$.}} By the Chernoff bound (b) with $\delta=1/2$,
\[
 \P\Bigl(N_n>\frac{\mu_n}2\Bigr)\geq 1-e^{-\frac{\mu_n}8}\geq 1-e^{-\frac 14}.
\]
For $m>\mu_n/2$ we have $m\geq 2$, since $\mu_n> 2$; hence $m(m-1)\geq m^2/2>\mu_n^2/8$, and Lemma~\ref{lem:birthday} gives
\[
 q_{n,m}=\P(E_n\mid N_n=m)\geq c_0\min\Bigl\{1,\frac{\mu_n^2}{8n}\Bigr\}\geq\frac{c_0}8\min\Bigl\{1,\frac{\mu_n^2}n\Bigr\}.
\]
Therefore, by \eqref{eq:birthday} and the monotonicity of $q_{n,m}$
\begin{align*}
 \P(E_n)&\geq\sum_{m>\mu_n/2} \P(N_n=m)\, q_{n,m}\geq \P\Bigl(N_n>\frac{\mu_n}2\Bigr) \, \frac{c_0}8\min\Bigl\{1,\frac{\mu_n^2}n\Bigr\}\\
 &\geq\frac{c_0\bigl(1-e^{-1/4}\bigr)}8\min\Bigl\{1,\frac{\mu_n^2}n\Bigr\}.
\end{align*}

Both cases give \eqref{eq:lower} with $c=\min\bigl\{e^{-2}/4,\ c_0(1-e^{-1/4})/8\bigr\}$.

Notice that the proof of  Case 2 does not use the estimates of Lemma~\ref{lem:pnk}.
\end{proof}

\section*{Appendix. The case where $\rho(x)$ does not tend to $\infty$ as $x\to\infty$.}%\label{ann:rho-infty}
Let us see here why if condition \eqref{eq:rho-to-inft} does not hold then
$\P(\Lambda_\phi^M\ \text{separated})=0$.
Assume there exist $C>0$ and $x_k\in\R$, $k\geq 1$, such that $\lim\limits_{k\to\infty} x_k=+\infty$ and $\rho(x_k)\leq C$. 

For any fixed $l\in\N$ consider the thin annuli around $|z|=x_k$
\[
 A_k^l=\Bigl\{z\in\C : x_k-\frac 1l\rho(x_k)<|z|\leq x_k+\frac 1l\rho(x_k)\Bigr\}.
\]
Taking a subsequence of $(x_k)_k$ if necessary, we can assume that the $A_k^l$ are pairwise disjoint. Notice that, by Lemma~\ref{lem:rho}(c),
\begin{equation}\label{eq:mukl}
 \mu_{k,l}:=\E[N(A_k^l)]\simeq\int_{A_k^l}\frac{dm(z)}{\rho^2(z)}\simeq \frac{|A_k^l|}{\rho^2(x_k)}\simeq\frac{\rho(x_k)\, x_k}{l \rho^2(x_k)}=\frac{x_k}{l\rho(x_k)}.
\end{equation}
By Lemma~\ref{lem:rho}(b), this value goes to infinity as $k\to \infty$ faster than a certain power of $x_k$.

Now split $A_k^l$ in angular sectors $Q_j^{k,l}$ of equal size comparable to $\rho(x_k)/l$: take $n_{k,l}=\bigl[ l x_k/\rho(x_k) \bigr]$ sectors of equal angle $(2\pi)/n_{k,l}$ and let us see that the condition 
\begin{equation}\label{eq:final}
 \frac{\mu_{k,l}^2}{n_{k,l}}\simeq\frac 1{l^3}\frac{x_k}{\rho(x_k)}\ \stackrel{k\to\infty}{\longrightarrow}\ \infty
\end{equation}
implies that the collision probability of having two points in the same sector is bounded below by a constant. 

Let $X_j^{k,l}=N(Q_j^{k,l})$ and consider the collision events
\[
E_{k,l}=\bigl\{\exists\, j=1,\dots, n_{k,l}\; :\; X_j^{k,l}\geq 2\bigr\}=\bigl\{\max_j X_{j}^{k,l}\geq 2\bigr\}.
\]
Then, for $\delta>0$ small, and denoting $N_{k,l}=N(A_{k,l})$,
\begin{align*}
 \P(E_{k,l})&=\sum_{m=0}^\infty \P(E_{k,l}\mid N_{k,l}=m)\P(N_{k,l}=m)\\
 &\geq \sum\limits_{m: |m-\mu_{k,l}|<\delta \mu_{k,l}} \P(E_{k,l}\mid N_{k,l}=m)\, \P(N_{k,l}=m)\\
 &\geq \sum\limits_{m: |m-\mu_{k,l}|<\delta \mu_{k,l}}\, \P\bigl(E_{k,l}\mid N_{k,l}=[(1-\delta)\mu_{k,l}]+1\bigr)\, \P(N_{k,l}=m)\\
 &=\P(E_{k,l}\mid N_{k,l}=[(1-\delta)\mu_{k,l}]+1) \, \P\bigl(|N_{k,l}-\mu_{k,l}|<\delta \mu_{k,l}\bigr)
\end{align*}

Since $\mu_{k,l}$ grows to infinity as $k$ increases, we can use Chernoff's estimate (c) in the second factor to see that for $k$ large
\begin{equation}\label{eq:app}
 \P(E_{k,l})\geq \frac 12\, \P(E_{k,l}\mid N_{k,l}=[(1-\delta)\mu_{k,l}]+1)
\end{equation}
From here we can proceed as in Section~\ref{sec:3.2}: letting $m=[(1-\delta)\mu_{k,l}]+1$ and $n=n_{k,l}$, the birthday probabilities 
\eqref{eq:birthday} yield, for $m\leq n$,
\begin{align*}
 1- \P(E_{k,l}\mid N_{k,l}=m)&=\prod_{j=1}^{m-1}\bigl(1-\frac jn\bigr)\leq \exp\bigl(-\frac{m(m-1)}{2n}\bigr).
\end{align*}
Since, by \eqref{eq:final},
\[
 \frac{m(m-1)}n=\frac{([(1-\delta)\mu_{k,l}]+1)\, [(1-\delta)\mu_{k,l}]}{n_{k,l}}\simeq \frac{\mu_{k,l}^2}{n_{k,l}}
\]
tends to infinity as $k\to\infty$, we deduce from \eqref{eq:app} that 
\[
\liminf_{k\to\infty}\P(E_{k,l})\geq \frac 12 .
\]
All combined, $\sum_{k=1}^\infty \P(E_{k,l})=+\infty$ for any $l\in\N$.
Following exactly the same procedure as in Proposition~\ref{prop:na} we see also that the events $\{E_{k,l}\}_k$ are negatively correlated, for any $l$ fixed. Hence, by Lemma~\ref{lem:BC2}, $\P(E_{k,l}\ \text{i.o.})=1$, for any $l$, and $\Lambda_\phi^M$ is almost surely not separated.

\subsection*{Use of AI tools} An earlier version of this manuscript was reviewed with the assistance of ChatGPT-6 (OpenAI) and Claude Opus 5.5 (Anthropic), which helped identify a correctable error and provided suggestions for improving the organization and clarity of the presentation. The authors have retained the original mathematical approach and overall structure of the manuscript and take full responsibility for its contents.

\bibliographystyle{plain}

\end{document}